\documentclass[11pt]{amsart}
\usepackage[T1]{fontenc}
\usepackage{lmodern}
\usepackage{microtype}
\usepackage[a4paper,margin=29mm]{geometry}
\usepackage{amssymb,amsmath}
\usepackage[hidelinks]{hyperref}
\usepackage{hyperref}
\hypersetup{pdftitle={Non-closed products in branch and weakly branch groups},
pdfauthor={Daniele D'Angeli and Emanuele Rodaro}}

\newtheorem{theorem}{Theorem}[section]
\newtheorem{lemma}[theorem]{Lemma}
\newtheorem{proposition}[theorem]{Proposition}
\newtheorem{corollary}[theorem]{Corollary}
\newtheorem{thmx}{Theorem}

\theoremstyle{definition}
\newtheorem{definition}[theorem]{Definition}
\theoremstyle{remark}
\newtheorem{remark}[theorem]{Remark}
\newtheorem{question}[theorem]{Question}

\DeclareMathOperator{\Aut}{Aut}
\DeclareMathOperator{\St}{St}
\DeclareMathOperator{\rist}{rist}
\DeclareMathOperator{\Fin}{Fin}
\DeclareMathOperator{\Fix}{Fix}
\DeclareMathOperator{\End}{End}
\newcommand{\cl}[1]{\overline{#1}}
\newcommand{\T}{\mathcal{T}}
\newcommand{\Lv}{\mathcal{L}}
\newcommand{\Gr}{\mathcal{G}}
\newcommand{\B}{\mathcal{B}}

\begin{document}

\title[Non-closed products in branch groups]{Non-closed products in branch\\
and weakly branch groups}
\author{Daniele D'Angeli}
\address{Universit\`a Niccol\`o Cusano, Rome, Italy}
\email{daniele.dangeli@unicusano.it}
\author{Emanuele Rodaro}
\address{Politecnico di Milano, Milan, Italy}
\email{emanuele.rodaro@polimi.it}

\subjclass[2020]{Primary 20E26, 20E08; Secondary 20E18, 20E07, 20E45, 20F10}
\keywords{(Weakly) branch groups, profinite topology,
product separability, congruence subgroup property,
activity}

\begin{abstract}
We prove that the first Grigorchuk group and the Gupta--Sidki $3$-group are not
double coset separable, disproving a conjecture of Grigorchuk, Leemann and
Nagnibeda. More generally, every finitely generated regular branch group with
bounded activity contains two finitely generated, profinitely closed subgroups
with trivial intersection whose product is not closed. Every infinite finitely
generated subgroup of the first Grigorchuk group or the Gupta--Sidki $3$-group
also fails to be double coset separable. We extend the construction to the
$p$-Basilica groups and, using polynomial activity, to the long-range group,
whose derived subgroup we prove to be finitely generated. Under an additional
hypothesis on finite quotients, we determine the closure of the product.
Membership in the constructed product is decidable in the first Grigorchuk group.
\end{abstract}

\maketitle
\raggedbottom

\section{Introduction}\label{sec:intro}

Let $G$ be a residually finite group. The \emph{profinite topology} on $G$ is the
group topology having the cosets of the normal subgroups of finite index as a
basis (see e.g.\ \cite{RZ00}); in other words, a subset $Y\subseteq G$ is closed if
every element outside $Y$ can be separated from $Y$ in a finite quotient of $G$.
Recall that $G$ is \emph{subgroup separable} if every finitely generated
subgroup of $G$ is closed, and that $G$ is \emph{double coset separable} if the
product $H_1H_2$ is closed for all finitely generated subgroups $H_1,H_2\le G$;
since translations are homeomorphisms, this is equivalent to requiring that every
double coset $H_1gH_2$ of finitely generated subgroups be closed. More generally,
$G$ is \emph{product separable} if every finite product of finitely generated
subgroups of $G$ is closed.

Product separability implies double coset separability,
which in turn implies subgroup separability (take $H_2=1$). Neither converse
holds. Polycyclic groups are double coset separable, but a product of three
finitely generated subgroups of the integral Heisenberg group need not be closed
\cite{LW80}; moreover, Minasyan constructed subgroup separable groups which are not
double coset separable \cite{Min23}. Free groups, on the other hand, are double coset
separable \cite{GR95,Nib92} and even product separable \cite{RZ93}; more recently,
limit groups have been shown to be product separable, and finitely generated
Kleinian groups to be double coset separable \cite{MM25}.

The first Grigorchuk group $\Gr$ and the Gupta--Sidki $3$-group $\Gr_3$ are
subgroup separable \cite{GW03,Gar16}, and so is every finitely generated
self-similar branch group with the subgroup induction property \cite{FL25}.
The finitely generated subgroups of $\Gr$ and $\Gr_3$ are, up to finite index,
products of diagonal copies of finite-index subgroups of the ambient group
\cite{GLN21,FGLN25}. In \cite[Section~2.1]{GLN21}, Grigorchuk, Leemann and Nagnibeda
conjectured that $\Gr$ and $\Gr_3$ are product separable.

We disprove this conjecture. In fact, it fails already for products of two
subgroups, so that $\Gr$ and $\Gr_3$ are not even double coset separable, and
this failure occurs throughout the class of regular branch groups with bounded
activity. Non-closed subgroups of weakly branch groups have also been studied in
\cite{FLN25}; those examples are infinitely generated, whereas here both factors
are finitely generated and closed.

In order to state our results, let us fix some notation; we refer
to Section~\ref{sec:prelim} for the details. Let $T_d=X^*$ be the regular rooted tree
over $X=\{0,\ldots,d-1\}$, where $d\ge2$. For a vertex $v$ we denote by
$\lambda_v(g)$ the automorphism acting as $g$ on the subtree rooted at $v$ and
fixing every vertex outside it, and we put $\lambda=\lambda_0$ and
$\rho=\lambda_{d-1}$; powers of $\lambda$ and $\rho$ denote compositions, so that,
for instance, $\rho^2=\lambda_{(d-1)(d-1)}$. The \emph{section} of an automorphism
at a vertex is its induced action between the corresponding rooted subtrees,
identified with $T_d$. An action has \emph{bounded activity} if each element has a
bounded number of non-trivial sections at each level. An automorphism is
\emph{finitary} if all its sections at some level are trivial. We write $K'$ for
the derived subgroup of $K$, and all closures are taken in the profinite topology
of the indicated group.

\begin{thmx}\label{thm:A}
Let $G\le\Aut(T_d)$ be finitely generated and level-transitive, with bounded
activity, and suppose that $G$ is regular branch over $K$
(Definition~\ref{def:branch}). Put
\[
 P=\{\lambda(x)\rho^2(x):x\in K\},\qquad
 Q=\{\lambda(y)\rho(y):y\in K\},\qquad
 W=\lambda(K)\rho(K).
\]
Then $P\cong Q\cong K$ are finitely generated and closed in $G$, and
$P\cap Q=1$. The subgroup $W\cong K\times K$ is closed in $G$, and it has finite
index in $G$ if and only if $d=2$. Moreover, $K'$ contains a non-finitary element,
and every such element $k$ satisfies
\[
 t_k:=\lambda^2(k)\in\cl{PQ}\setminus PQ.
\]
In particular, $G$ is not double coset separable.
\end{thmx}

Notice that $Q=\{\lambda_0(y)\lambda_{d-1}(y):y\in K\}$ and
$P=\{\lambda_0(x)\lambda_{(d-1)(d-1)}(x):x\in K\}$ are diagonal copies of $K$
supported on the pairs of vertices $\{0,d-1\}$ and $\{0,(d-1)(d-1)\}$; in the
terminology of \cite{GLN21,FGLN25}, they are among the simplest diagonal block
subgroups. Recall that $G$ has the \emph{congruence subgroup property} if every
subgroup of finite index of $G$ contains the pointwise stabiliser of some level,
and that $G$ has \emph{trivial branch kernel} if every normal subgroup of finite
index of $G$ contains the rigid stabiliser of some level
(Section~\ref{sec:prelim}); the former property implies the latter. Under this
additional assumption the restriction $k\in K'$ can be removed, and we can
determine the closure of the product.

\begin{thmx}\label{thm:B}
If $G$ in Theorem~\ref{thm:A} has trivial branch kernel (for instance, if $G$ has
the congruence subgroup property), then $\cl{PQ}=W$. Furthermore,
$t_k\in\cl{PQ}\setminus PQ$ for every non-finitary $k\in K$.
\end{thmx}

The construction also works under weaker branching assumptions. The next theorem
separates the non-closedness of the product from the additional conditions needed
to make its two factors closed. Here $G$ is \emph{self-similar} if every section of
an element of $G$ belongs again to $G$.

\begin{thmx}\label{thm:C}
Let $G\le\Aut(T_d)$ be level-transitive with bounded activity. Suppose that
$1\ne K\le G$ is finitely generated and that $\lambda_x(K)\le K$ for every
$x\in X$. Define $P$ and $Q$ as in Theorem~\ref{thm:A}. Then $P\cong Q\cong K$
and $P\cap Q=1$. Moreover, $K'$ contains a non-finitary element, and every such
element $k$ satisfies
\[
 \lambda^2(k)\in\cl{PQ}\setminus PQ.
\]
Hence $G$ is not double coset separable. If, in addition, $G$ is self-similar and
$K$ is closed in $G$, then $P$ and $Q$ are closed in $G$.
\end{thmx}

Theorem~\ref{thm:C} applies to the Basilica group and, more generally, to every
$p$-Basilica group, whose derived subgroup is finitely generated. For the
Basilica group this was proved by Francoeur
\cite[Proposition~4.7]{Fra20}; in Section~\ref{sec:basilica} we give an argument
which works for every prime $p$ (Lemma~\ref{lem:pbasilica-derived}), and we deduce
that neither the Basilica group nor any $p$-Basilica group is double coset
separable (Corollary~\ref{cor:basilica}). For the groups appearing in the
conjecture of Grigorchuk, Leemann and Nagnibeda we obtain the following.

\begin{corollary}\label{cor:GLN}
The first Grigorchuk group and the Gupta--Sidki $3$-group are subgroup separable
but not double coset separable. The conjecture of \cite[Section~2.1]{GLN21}
therefore fails. Moreover, every infinite finitely generated subgroup of either
group is subgroup separable but not double coset separable.
\end{corollary}

The last assertion follows from the description of the finitely generated
subgroups of these groups obtained in \cite{FGLN25}, and holds more generally for
every torsion GGS-group (Corollary~\ref{cor:hered}).

Bounded activity is used only to show that a certain recursively defined
automorphism does not belong to $G$, and it can be replaced by polynomial activity,
provided that a suitable element of $K'$ has activity of the same degree as $G$
(Theorem~\ref{thm:Cpoly}). We illustrate this with a weakly branch group of linear
activity, namely the long-range group, whose Schreier graphs were studied by
Benjamini and Hoffman \cite{BH05} and whose amenability was proved by Amir, Angel
and Vir\'ag \cite{AAV13}. In this example the finite generation of the branching
subgroup has to be established directly, and the witness has to be chosen with
linear activity.

\begin{thmx}\label{thm:D}
Let $G=\langle a,b\rangle\le\Aut(T_2)$, where $a=(1,a)(0\ 1)$ and $b=(a,b)$. Then
$G$ is isomorphic to the long-range group and has linear activity. It is weakly
regular branch over its finitely generated derived subgroup $G'$, but this action
is not branch. The subgroups
\[
 P=\{\lambda(x)\rho^2(x):x\in G'\},\qquad Q=\{\lambda(y)\rho(y):y\in G'\}
\]
are finitely generated, closed and isomorphic to $G'$, with $P\cap Q=1$, while
$PQ$ is not closed.
In particular, the long-range group is not double coset separable.
\end{thmx}

Non-finitarity alone no longer gives an activity obstruction: the generator $b$
itself solves $x=\lambda(a)\rho(x)$, with $a$ non-finitary
(Remark~\ref{rem:linear-witness}).

Let us briefly describe the strategy of the proofs. The non-membership of $t_k$ in
$PQ$ is reduced to the recursive equation $x=\lambda(k)\rho(x)$, whose
unique solution in $\Aut(T_d)$ has unbounded activity when $k$ is non-finitary
(Lemma~\ref{lem:activity}); this is the only point where bounded activity is used,
and in Section~\ref{sec:linear} we replace it by a comparison of activity growth.
On the other hand, finite truncations of this solution approximate $t_k$ by elements
of $PQ$ in every finite quotient of $G$. This relies on a classical observation of
Grigorchuk: every non-trivial normal subgroup of a level-transitive group contains
the derived subgroups of the rigid stabilisers of the vertices of some level
(Lemma~\ref{lem:grig}). This is the reason why the witness $k$ is taken in $K'$, and
in Section~\ref{sec:hanoi} we show, using the Hanoi towers group, that this
restriction cannot be dropped in general. From a different point of view, the
product $PQ$ corresponds to the product of the graph of $\rho$ and of the diagonal
subgroup of $K\times K$, so that Theorem~\ref{thm:A} says that the $\rho$-twisted
conjugacy class of the identity of $K$ is not closed (Section~\ref{sec:twisted}).

After the preliminaries, Section~\ref{sec:proofs} proves
Theorems~\ref{thm:A}--\ref{thm:C}. Sections~\ref{sec:examples}
and~\ref{sec:basilica} give the applications to regular branch groups and the
Basilica groups. Section~\ref{sec:linear} extends the argument to polynomial
activity and proves Theorem~\ref{thm:D}. The Hanoi towers example in
Section~\ref{sec:hanoi} explains the restriction on the witness in
Theorem~\ref{thm:A}. Section~\ref{sec:twisted} treats twisted conjugacy, rational
subsets and the membership algorithm. We conclude with questions in
Section~\ref{sec:activity}.

\section{Preliminaries}\label{sec:prelim}

We recall the definitions and facts used in the proofs. For further background
on branch groups, see \cite{BGS03}.

\subsection{Sections and rigid embeddings}
The vertices of $T_d=X^*$ are the finite words over $X$, the root being the empty
word, and there is an edge from $v$ to $vx$ for every $x\in X$. We denote by
$\Lv_n=X^n$ the $n$-th level of $T_d$ and by $\T_v=vX^*$ the subtree rooted at $v$.
Two vertices are \emph{incomparable} if neither of them is a prefix of the other. We
use right actions, so that $v^{gh}=(v^g)^h$; accordingly, we put $g^h=h^{-1}gh$ and
$[g,h]=g^{-1}h^{-1}gh$, so that $G'=\langle[g,h]:g,h\in G\rangle$. For
$g\in\Aut(T_d)$ and a vertex $v$, the \emph{section} $g|_v\in\Aut(T_d)$ is defined by
$(vu)^g=v^g u^{g|_v}$ for every $u\in X^*$, and one has
\begin{equation}\label{eq:sections}
 (gh)|_v=g|_v\,h|_{v^g},\qquad
 (g^{-1})|_{v^g}=(g|_v)^{-1}.
\end{equation}
We write $g=(g|_0,\ldots,g|_{d-1})\sigma$ if $\sigma$ is the permutation induced by
$g$ on $X$; a \emph{rooted permutation} is an automorphism whose sections at the
first level are all trivial. An action is \emph{level-transitive} if it is
transitive on every level $\Lv_n$, and $G\le\Aut(T_d)$ is \emph{self-similar} if
$g|_v\in G$ for every $g\in G$ and every vertex $v$.

For $G\le\Aut(T_d)$ we denote by $\St_G(v)$ the stabiliser of the vertex $v$ and by
$\St_G(n)$ the pointwise stabiliser of $\Lv_n$; we abbreviate
$\St_{\Aut(T_d)}(n)$ to $\St(n)$. The map
\[
 \psi\colon\St(1)\longrightarrow\Aut(T_d)^d,
 \qquad g\longmapsto(g|_0,\ldots,g|_{d-1}),
\]
is an isomorphism. The \emph{rigid embedding} $\lambda_v$ is defined by
$(vu)^{\lambda_v(g)}=vu^g$ for $u\in X^*$, while $\lambda_v(g)$ fixes every vertex
outside $\T_v$. It is an injective homomorphism, with
$\lambda_{vw}=\lambda_v\circ\lambda_w$ and $\lambda_v(g)\in\St(|v|)$; moreover, the
elements of $\lambda_v(\Aut(T_d))$ and of $\lambda_w(\Aut(T_d))$ commute whenever
$v$ and $w$ are incomparable. The \emph{rigid stabiliser} of $v$ in $G$ is
$\rist_G(v)=G\cap\lambda_v(\Aut(T_d))$, and the rigid stabiliser of the $n$-th level
is $\rist_G(n)=\prod_{v\in\Lv_n}\rist_G(v)$. The support of an automorphism is
the set of vertices it moves.

\subsection{Activity}
For $g\in\Aut(T_d)$ we put $A_n(g)=\#\{v\in\Lv_n:g|_v\ne1\}$. By
\eqref{eq:sections} we have
\begin{equation}\label{eq:subadd}
 A_n(gh)\le A_n(g)+A_n(h),\qquad A_n(g^{-1})=A_n(g).
\end{equation}
The group $G$ has \emph{bounded activity} if $\sup_n A_n(g)<\infty$ for every
$g\in G$; by \eqref{eq:subadd}, for a finitely generated group it is enough to check
this on a finite set of generators. More generally, the action has
\emph{polynomial activity of degree at most $m$} if $A_n(g)=O(n^m)$ for every
$g\in G$, with constants allowed to depend on $g$. We use \emph{linear activity}
for the case of degree at most $1$, but not bounded activity.
An automorphism $g$ is \emph{finitary} if
$A_n(g)=0$ for some $n$. The groups $\Fin_n=\{g:A_n(g)=0\}$ are finite and form an
increasing sequence, so that their union $\Fin$ is locally finite, that is, each of
its finitely generated subgroups is finite. An automorphism is
\emph{finite-state} if it has only finitely many distinct sections. In this case it
can be represented by a finite graph (a Mealy automaton) whose states are its
sections, each state being labelled by its root permutation and having transitions
to its first-level sections.

\begin{lemma}\label{lem:nonfin}
Every infinite finitely generated subgroup of $\Aut(T_d)$ contains a non-finitary
element. If $k$ is non-finitary, then $A_n(k)\ge1$ for every $n$.
\end{lemma}

\begin{proof}
A finitely generated subgroup consisting of finitary elements is contained in some
$\Fin_n$, hence it is finite. The second assertion is just the definition of a
non-finitary automorphism.
\end{proof}

\subsection{Branching}
All branching properties below refer to the specified action on the rooted tree.
Recall that a level-transitive group $G\le\Aut(T_d)$ is \emph{weakly branch} if
$\rist_G(v)\ne1$ for every vertex $v$, and it is \emph{branch} if $\rist_G(n)$ has
finite index in $G$ for every $n$. We use the following version of regular
branching, which does not require self-similarity and allows the branching subgroup
not to fix the first level.

\begin{definition}\label{def:branch}
Let $G\le\Aut(T_d)$ be level-transitive, and let $1\ne K\le G$. The group $G$ is
\emph{weakly regular branch over $K$} if
\[
 K^d\le\psi\bigl(K\cap\St_G(1)\bigr).
\]
It is \emph{regular branch over $K$} if, in addition, $K$ is normal of finite
index in $G$ and $[\psi(K\cap\St_G(1)):K^d]<\infty$.
\end{definition}

Notice that the inclusion $K^d\le\psi(K\cap\St_G(1))$ is equivalent to
$\lambda_x(K)\le K$ for every $x\in X$, and it implies $\lambda_v(K)\le K$ for every
vertex $v$. In the regular branch case, the subgroup
\[
 W':=\psi^{-1}(K^d)=\prod_{x\in X}\lambda_x(K)\cong K^d
\]
has finite index in $K$, since
\[
 [K:W']=[K:K\cap\St_G(1)]\,[\psi(K\cap\St_G(1)):K^d];
\]
iterating this observation, one sees that $G$ is a branch group. If $G$ is
self-similar and $[G:K]<\infty$, the second index in the definition is
automatically finite, since $K^d$ has finite index in $G^d$. These conventions agree
with the usual ones in the self-similar setting \cite{BGS03,BSZ12}.

The following lemma, due to Grigorchuk, allows us to work with arbitrary finite
quotients, without assuming the congruence subgroup property. For the reader's
convenience we include its short proof.

\begin{lemma}[{\cite[Theorem~4]{Gri00}, \cite[Lemma~2.2]{BSZ12}}]\label{lem:grig}
Let $G\le\Aut(T_d)$ be level-transitive and let $1\ne N\trianglelefteq G$. Then
there exists $m$ such that $\rist_G(u)'\le N$ for every $u\in\Lv_m$.
\end{lemma}

\begin{proof}
Let us choose $1\ne g\in N$ and a vertex $w$ with $w^g\ne w$. If
$x,y\in\rist_G(w)$, then $(x^{-1})^g$ is supported in $\T_{w^g}$, so that it commutes
with $x$ and $y$. Hence
\[
 [[g,x],y]=[(x^{-1})^gx,y]=[x,y]\in N,
\]
that is, $\rist_G(w)'\le N$. Since $\rist_G(w)^h=\rist_G(w^h)$ for every $h\in G$,
level-transitivity gives the assertion for every vertex of level $m=|w|$.
\end{proof}

We say that $G$ has the \emph{congruence subgroup property} if every subgroup of
finite index of $G$ contains $\St_G(n)$ for some $n$, and that $G$ has
\emph{trivial branch kernel} if every normal subgroup of finite index of $G$
contains $\rist_G(n)$ for some $n$ \cite{BSZ12}. Since $\rist_G(n)\le\St_G(n)$, the
congruence subgroup property implies that the branch kernel is trivial.

\subsection{Profinite topology}
We use the standard facts about the profinite topology recorded in \cite{RZ00}.
Every subgroup of $\Aut(T_d)$ is residually finite, since the intersection of its
level stabilisers is trivial. A subgroup of finite index is open and closed, and its
own profinite topology coincides with the induced one; for an arbitrary subgroup,
instead, its own profinite topology may be finer than the induced one. Every
homomorphism is continuous with respect to the profinite topologies, and the
profinite topology of a finite direct product is the product of the profinite
topologies of its factors.

The topology on $G$ defined by the level stabilisers $\St_G(n)$ is called the
\emph{congruence topology}. It is coarser than the profinite topology of $G$, and
the two topologies coincide precisely when $G$ has the congruence subgroup property.
In particular, a subset of $G$ which is closed in the congruence topology is closed
in the profinite topology.

\begin{lemma}\label{lem:graph}
Let $\alpha\colon H_1\to H_2$ be a homomorphism, where $H_2$ is residually finite.
Then the graph $\Gamma_\alpha=\{(x,\alpha(x)):x\in H_1\}$ is closed in $H_1\times H_2$.
\end{lemma}

\begin{proof}
Let $(u,v)\notin\Gamma_\alpha$, and let $\theta\colon H_2\to F$ be a homomorphism
to a finite group with $\theta(v)\ne\theta(\alpha(u))$. The inverse image of the
diagonal of $F\times F$ under $(x,y)\mapsto(\theta\alpha(x),\theta(y))$ is a
subgroup of finite index of $H_1\times H_2$ containing $\Gamma_\alpha$ but not
$(u,v)$.
\end{proof}

The next observation will provide closed factors in the case of the Basilica
groups.

\begin{lemma}\label{lem:rigid-diagonal}
Let $G\le\Aut(T_d)$ be self-similar and let $K\le G$ be closed. Suppose that $V$ is
a non-empty finite set of pairwise incomparable vertices such that
$\lambda_v(K)\le G$ for every $v\in V$. Then the following subgroup is closed in $G$:
\[
 \Bigl\{\prod_{v\in V}\lambda_v(k):k\in K\Bigr\}.
\]
\end{lemma}

\begin{proof}
Let us consider the subgroup of finite index $U=\bigcap_{v\in V}\St_G(v)$. By
self-similarity, for every $v\in V$ the section map $g\mapsto g|_v$ is a
homomorphism $U\to G$, which is therefore continuous. The set of elements of $U$
whose sections at the vertices of $V$ coincide is closed, by
Lemma~\ref{lem:graph} applied to the identity of $G$, and the set of those whose
section at a fixed $v\in V$ belongs to $K$ is closed, since $K$ is. Finally, the set
of elements of $U$ fixing every vertex outside $\bigcup_{v\in V}\T_v$ is an
intersection of vertex stabilisers, each of which has finite index, hence it is
closed. The intersection of these three closed sets is exactly the displayed
subgroup, which is therefore closed in $U$, and hence in $G$.
\end{proof}

Finally, let us record how double coset separability behaves with respect to
subgroups and finite-index overgroups.

\begin{lemma}\label{lem:hered}
Subgroups and finite-index overgroups of double coset separable groups are double
coset separable.
\end{lemma}

\begin{proof}
Let $G$ be double coset separable, let $L\le G$, and let $H_1,H_2\le L$ be finitely
generated. Then $H_1H_2$ is closed in $G$, hence in the topology induced on $L$,
and a fortiori in the profinite topology of $L$.

Conversely, suppose that $L\le G$ has finite index and is double coset separable,
and let $H_1,H_2\le G$ be finitely generated. The subgroups $H_i\cap L$ have finite
index in $H_i$, hence they are finitely generated, and we can choose finite coset
decompositions $H_1=\bigcup_i h_i(H_1\cap L)$ and $H_2=\bigcup_j(H_2\cap L)k_j$.
Then
\[
 H_1H_2=\bigcup_{i,j}h_i\bigl((H_1\cap L)(H_2\cap L)\bigr)k_j
\]
is a finite union of translates of subsets which are closed in $L$, hence in $G$.
Therefore $H_1H_2$ is closed in $G$.
\end{proof}

\section{Proofs of Theorems A--C}\label{sec:proofs}

The first lemma gives the obstruction to membership in $PQ$.

\begin{lemma}\label{lem:activity}
Let $k\in\Aut(T_d)$ be non-finitary. If $x\in\Aut(T_d)$ satisfies
$x=\lambda(k)\rho(x)$, then
\[
 A_n(x)=1+\sum_{j=0}^{n-1}A_j(k)\ge n+1\qquad(n\ge1).
\]
\end{lemma}

\begin{proof}
The equation says that $\psi(x)=(k,1,\ldots,1,x)$. Consequently,
$A_n(x)=A_{n-1}(k)+A_{n-1}(x)$ for every $n\ge1$. Since $k\ne1$, we also have
$x\ne1$, that is, $A_0(x)=1$. Iterating, we obtain the formula, and the inequality
follows from Lemma~\ref{lem:nonfin}.
\end{proof}

Let us note that the equation $x=\lambda(k)\rho(x)$ has a unique solution in
$\Aut(T_d)$, namely the infinite product $\prod_{j\ge0}\lambda_{(d-1)^j0}(k)$. Indeed,
its factors have pairwise disjoint supports and fix successively deeper levels, so
that the product converges in the congruence topology, and the recursion determines
the action on each finite level. If $k$ is finitary, the sum in
Lemma~\ref{lem:activity} is eventually constant; the assumption on $k$ is therefore
essential for the activity argument.

The second lemma guarantees the existence of suitable witnesses in the derived
subgroup.

\begin{lemma}\label{lem:Kprime}
Let $1\ne K\le\Aut(T_d)$ be finitely generated and suppose that $\lambda_v(K)\le K$
for every vertex $v$. Then $K$ is infinite and non-abelian, and $K'$ contains a
non-finitary element.
\end{lemma}

\begin{proof}
Since $K$ contains $\prod_{v\in\Lv_n}\lambda_v(K)\cong K^{d^n}$ for every $n$, it
cannot be finite. Nor can it be abelian: being finitely generated, it would be
isomorphic to $\mathbb Z^r\times F$ with $F$ finite, and an embedding of $K^d$ into
$K$ would force both $r$ and $F$ to be trivial. Thus $K'\ne1$.

Let us choose a non-finitary $h\in K$, which exists by Lemma~\ref{lem:nonfin}, and
an element $1\ne x\in K'$. If $x$ is non-finitary there is nothing to prove.
Otherwise, take $m$ with $x\in\Fin_m$, choose $p\in\Lv_m$ moved by $x$, and put
$q=p^{x^{-1}}\ne p$ and $g=\lambda_p(h)\in K$. The automorphism $g^{-1}$ fixes
$q$ with trivial section, while $q^x=p$ and $x|_q=1$. Thus
\eqref{eq:sections} gives $(x^g)|_q=h$. Consequently $x^g$ is non-finitary and
belongs to $K'$, since $K'$ is normal in $K$.
\end{proof}

\begin{proof}[Proof of Theorem~\ref{thm:A}]
Let $\iota(x,y)=\lambda(x)\rho(y)$ for $x,y\in K$. Since $\lambda(K)$ and $\rho(K)$
have disjoint supports, $\iota\colon K\times K\to W$ is an isomorphism; under this
isomorphism, $P$ corresponds to the graph of $\rho\colon K\to K$ and $Q$ to the
diagonal subgroup. In particular, both are isomorphic to $K$, which is finitely
generated because it has finite index in $G$.

In $K^d\cong W'$, the inverse images of $P$ and
$Q$ are, respectively,
\[
 \{(x,1,\ldots,1,\rho(x)):x\in K\},\qquad
 \{(y,1,\ldots,1,y):y\in K\}.
\]
They are graphs of homomorphisms $K\to K^{d-1}$, hence they are closed by
Lemma~\ref{lem:graph}. The subgroup corresponding to $W$ is
$K\times1^{d-2}\times K$, which is closed as well, being the kernel of the
projection onto the residually finite group $K^{d-2}$. As $W'$ has finite index in
$G$, all three subgroups are closed in $G$. If $d=2$, then $W=W'$; if $d>2$, then
$[W':W]=|K|^{d-2}=\infty$ by Lemma~\ref{lem:Kprime}. Finally, an element of $P\cap Q$
satisfies $\iota(x,\rho(x))=\iota(y,y)$, hence $x=y=\rho(x)$; iterating, we obtain
$x=\rho^n(x)\in\St(n)$ for every $n$, so that $x=1$.

To prove that the product is not closed, we first show that a $t_k$ is not in $PQ$: qchoose a non-finitary
$k\in K'$, which exists by Lemma~\ref{lem:Kprime}, and put $c=\lambda(k)$, so that
$t_k=\iota(c,1)$. If $t_k\in PQ$, by injectivity of $\iota$, there are
$x,y\in K$ with $xy=c$ and $\rho(x)y=1$. Eliminating $y$ we obtain $x=c\rho(x)$,
which contradicts bounded activity by Lemma~\ref{lem:activity}. Notice that this
step works for every non-finitary $k\in K$.

To prove that $t_k$ is in the closure, we only use that $k\in K'$: recall that
$\lambda(K)\le K$ and $\rho(K)\le K$, so that $\lambda$ and $\rho$ restrict to
injective endomorphisms of $K$. In particular, $c=\lambda(k)\in K$ and
\[
 \rho^j(c)=\lambda_{(d-1)^j}\bigl(\lambda_0(k)\bigr)=\lambda_{(d-1)^j0}(k)\in K
 \qquad(j\ge0).
\]
Since the vertices $(d-1)^j0$, $j\ge0$, are pairwise incomparable, these
elements have pairwise disjoint supports, hence they commute. Putting
$x_n=\prod_{j=0}^{n-1}\rho^j(c)\in K$, we have
$\rho(x_n)=\prod_{j=1}^{n}\rho^j(c)$, and therefore
$x_n\rho(x_n)^{-1}=c\rho^n(c)^{-1}$. Since $x_n$ and $\rho(x_n)^{-1}$ belong to
$K$, the elements $\iota(x_n,\rho(x_n))$ and
$\iota(\rho(x_n)^{-1},\rho(x_n)^{-1})$ belong to $P$ and $Q$, respectively, so that
\begin{equation}\label{eq:approximation}
 \begin{split}
 &\iota(x_n,\rho(x_n))\,
    \iota(\rho(x_n)^{-1},\rho(x_n)^{-1})\\
 &\hspace{15mm}=\iota(c\rho^n(c)^{-1},1)
   =t_k\,\lambda(\rho^n(c))^{-1}\in PQ.
 \end{split}
\end{equation}
Let now $N\trianglelefteq G$ be a normal subgroup of finite index. Since $G$
is infinite, $N\ne1$, and Lemma~\ref{lem:grig} gives $m$ such that
$\rist_G(u)'\le N$ for every $u\in\Lv_m$. Let $n\ge m$, put
$v=0(d-1)^n0$ and let $u$ be the prefix of length $m$ of $v$. Since
$\lambda_v$ is a homomorphism, we have $\lambda_v(K')=\lambda_v(K)'$, and
since $\lambda_v(K)\le K\le G$ we have
$\lambda_v(K)\le\rist_G(v)\le\rist_G(u)$. As $k\in K'$, we obtain
\[
 \lambda(\rho^n(c))=\lambda_v(k)\in\lambda_v(K)'\le\rist_G(u)'\le N .
\]
By \eqref{eq:approximation}, $t_k\in PQ\,\lambda(\rho^n(c))\subseteq PQN$.
Since this holds for every normal subgroup $N$ of finite index, we conclude
that $t_k\in\cl{PQ}$. Notice that this step works for every $k\in K'$,
whereas the previous one works for every non-finitary $k\in K$; by
Lemma~\ref{lem:Kprime} there exist elements satisfying both conditions, and
for each of them we have $t_k\in\cl{PQ}\setminus PQ$.
\end{proof}

\begin{proof}[Proof of Theorem~\ref{thm:B}]
Let $N$ be a normal subgroup of finite index of $K$, and let
$N_0=\bigcap_{g\in G}N^g$ be its normal core in $G$. Since $K$ has finite index in
$G$, so does $N_0$. Since $G$ has
trivial branch kernel, there is $n\ge1$ with $\rist_G(n)\le N_0$, and hence
\[
 \rho^n(K)=\lambda_{(d-1)^n}(K)\le\rist_G\bigl((d-1)^n\bigr)\le\rist_G(n)\le N.
\]
For $z\in K$, put $x=z\rho(z)\cdots\rho^{n-1}(z)\in K$. Since $\rho$ is an
endomorphism of $K$, we have $\rho(x)=\rho(z)\rho^2(z)\cdots\rho^{n}(z)$, and
therefore $z\rho(x)=x\rho^n(z)$. Hence
\[
 x\rho(x)^{-1}=z\,\bigl(\rho(x)\rho^n(z)^{-1}\rho(x)^{-1}\bigr)\in zN,
\]
because $\rho^n(z)\in\rho^n(K)\le N$ and $N$ is normal in $K$. Since the cosets
$zN$, where $N$ ranges over the normal subgroups of finite index of $K$, form a
basis of neighbourhoods of $z$ in the profinite topology of $K$, this shows that
the set $D=\{x\rho(x)^{-1}:x\in K\}$ is dense in $K$.

The map $(u,v)\mapsto(uv^{-1},v)$ is a homeomorphism of $K\times K$, with
inverse $(z,v)\mapsto(zv,v)$. It sends the product of the graph of $\rho$ and
the diagonal onto $D\times K$, which is dense. Applying $\iota$, we see that
$PQ$ is dense in $W$ in its own profinite topology. The subgroup $W$ is a direct factor of $W'\cong K^d$, and
$W'$ has finite index in $G$. Thus the profinite topology of $W$ coincides with
the topology induced by $G$. Since $W$ is closed, we obtain
$\cl{PQ}=W$. Finally, the non-membership argument in the proof of
Theorem~\ref{thm:A} works for every non-finitary $k\in K$.
\end{proof}

\begin{proof}[Proof of Theorem~\ref{thm:C}]
In the proof of Theorem~\ref{thm:A}, the finite index of $K$ is used only to
ensure that $K$ is finitely generated and that $P$ and $Q$ are closed. We check
that the remaining arguments work.

By induction on $|v|$, the inclusions $\lambda_x(K)\le K$ give
$\lambda_v(K)\le K$ for every vertex $v$. Hence Lemma~\ref{lem:Kprime} provides a
non-finitary $k\in K'$, and $\lambda,\rho$ restrict to injective endomorphisms of
$K$. As before, $\iota(x,y)=\lambda(x)\rho(y)$ is an injective homomorphism
$K\times K\to G$, so that $P\cong Q\cong K$ are finitely generated, and
$P\cap Q=1$: if $\iota(x,\rho(x))=\iota(y,y)$, then $x=\rho(x)$, hence
$x=\rho^n(x)\in\St(n)$ for every $n$, that is, $x=1$.

The two halves of $\lambda^2(k)\in\cl{PQ}\setminus PQ$ use different hypotheses:
\begin{itemize}
\item $\lambda^2(k)\notin PQ$ only uses that $k$ is non-finitary and that $G$ has
bounded activity (Lemma~\ref{lem:activity});
\item $\lambda^2(k)\in\cl{PQ}$ only uses that $k\in K'$, that
$\lambda_v(K)\le\rist_G(v)$ for every $v$, and that $G$ is level-transitive
(Lemma~\ref{lem:grig}).
\end{itemize}
Thus both remain valid. Since $P$ and $Q$ are finitely generated and $PQ$ is not
closed, $G$ is not double coset separable.

Finally, if $G$ is self-similar and $K$ is closed, then
$P=\{\lambda_0(x)\lambda_{(d-1)(d-1)}(x):x\in K\}$ and
$Q=\{\lambda_0(y)\lambda_{d-1}(y):y\in K\}$ are closed by
Lemma~\ref{lem:rigid-diagonal}, applied with $V=\{0,(d-1)(d-1)\}$ and
$V=\{0,d-1\}$, respectively.
\end{proof}

\section{Examples}\label{sec:examples}

In this section we apply our results to several classes of regular branch groups.
We start with the first Grigorchuk group, for which Theorem~\ref{thm:B} gives an
explicit example.

\subsection{The first Grigorchuk group}
Let $\Gr=\langle a,b,c,d\rangle\le\Aut(T_2)$, where
\[
 a=(1,1)(0\ 1),\qquad b=(a,c),\qquad c=(a,d),\qquad d=(1,b).
\]
It is well known that $\Gr$ is regular branch over the normal closure
$K=\langle(ab)^2\rangle^{\Gr}$, which has index $16$ \cite[Chapter~VIII]{dlH00},
and that $\Gr$ has the congruence subgroup property \cite{Gri00}. The generators have
activity at most $2$. With our conventions we have $ab=(c,a)(0\ 1)$ and
$(ab)^2=(ca,ac)$; moreover, $(ab)^2|_{01}=d$, while $d|_1=b$, $b|_1=c$ and $c|_1=d$,
so that $(ab)^2$ is non-finitary. Therefore Theorem~\ref{thm:B} gives the explicit
example
\[
 \lambda^2((ab)^2)\in
 \cl{PQ}\setminus PQ
   =\psi^{-1}(K\times K)\setminus PQ.
\]
In particular, $PQ$ is dense in a subgroup of finite index of $\Gr$. Let us also
note that membership in this particular product is decidable
(Proposition~\ref{prop:membership}).

\subsection{Other regular branch groups}
Let us now describe the other families to which our results apply, in order to
specify the parameters involved. Let $p$ be an odd prime and let $a$ be the rooted
cycle $(0\ 1\ \cdots\ p-1)$. For a non-zero vector
$e=(e_1,\ldots,e_{p-1})\in\mathbb F_p^{p-1}$, the associated
\emph{Grigorchuk--Gupta--Sidki (GGS) group} is $\langle a,b\rangle$, where
$b=(a^{e_1},\ldots,a^{e_{p-1}},b)$. The Gupta--Sidki $3$-group corresponds to $p=3$
and $e=(1,-1)$. A \emph{multi-GGS group} is obtained by using several generators
$b_1,\ldots,b_r$ of this form, with linearly independent defining vectors.

If $e$ is non-symmetric, that is, $e_i\ne e_{p-i}$ for some $i$, and
$\sum_ie_i=0$, \emph{Pervova's group} associated with $e$ is $\langle a,b,c\rangle$,
where $b$ is as above and $c=(c,a^{e_1},\ldots,a^{e_{p-1}})$
\cite[Section~1.2]{Per07}.

For an odd integer $d\ge3$, let us identify $X$ with $\mathbb Z/d\mathbb Z$ and let
$\mu_i(x)=2i-x$. The \emph{Hanoihedral group} $D_d$ is generated by
$a_0,\ldots,a_{d-1}$, where $a_i$ has root permutation $\mu_i$, section $a_i$ at the
vertex $i$ and trivial sections at the other first-level vertices \cite{GS25}. For
$d=3$ this is, up to renaming the generators, the three-peg Hanoi towers group
discussed in Section~\ref{sec:hanoi}.

\begin{corollary}\label{cor:examples}
The following groups are not double coset separable:
\begin{enumerate}
\item GGS-groups with non-constant defining vector;
\item multi-GGS groups with at least two linearly independent defining vectors;
\item Pervova's groups with the parameters specified above;
\item the Hanoihedral groups $D_d$, for odd $d\ge3$;
\item the generalised Hanoi towers groups $G_n$, $n\ge3$, of Skipper
\cite{Ski20}, and their subgroups of finite index.
\end{enumerate}
In each case, the construction of Theorem~\ref{thm:A} gives two closed finitely
generated subgroups with trivial intersection and non-closed product.
\end{corollary}

\begin{proof}
All these groups are finitely generated, self-similar and level-transitive. A
GGS-group with non-constant defining vector is regular branch over
$\gamma_3(G)=[G',G]$ \cite{FAZR14}, and a multi-GGS group with at least two linearly
independent defining vectors is regular branch over its derived subgroup
\cite[Lemma~2]{GUA19}. Pervova's groups are regular branch over their derived
subgroups \cite[Proposition~1.3 and Lemma~1.4]{Per07}, and so are the Hanoihedral
groups \cite[Theorem~1]{GS25}. For the first three families, each non-rooted
generator has at most $p$ non-trivial sections on any level, while each Hanoihedral
generator has exactly one. Hence all these groups have bounded activity, and
Theorem~\ref{thm:A} applies. The GGS-groups and the multi-GGS groups in the
statement also have the congruence subgroup property \cite{FAGUA17,GUA19}, so that
Theorem~\ref{thm:B} applies to them as well. For $n\ge3$, Skipper's group $G_n=\langle a_1,\ldots,a_n\rangle\le\Aut(T_n)$, where
$a_i=(1,\ldots,1,a_i,1,\ldots,1)\sigma_i$ has the section $a_i$ at the vertex $i$ and
$\sigma_i$ is an $(n-1)$-cycle fixing $i$, is self-similar and level-transitive, and
it is regular branch over $G_n'$, which has finite index since
$G_n^{ab}\cong(\mathbb Z/(n-1)\mathbb Z)^n$ \cite{Ski20}; each of its generators has
exactly one non-trivial section on every level. Subgroups of finite index are covered
by Lemma~\ref{lem:hered}.
\end{proof}

\subsection{Infinite finitely generated subgroups}
Let $G\le\Aut(T_d)$ be \emph{self-replicating}, that is, self-similar with
$\{g|_v:g\in\St_G(v)\}=G$ for every vertex $v$. Following
\cite[Definitions~5.7 and~5.8]{FGLN25}, a collection $\mathcal X$ of subgroups of
$G$ is \emph{inductive} if \\
(I) $1,G\in\mathcal X$; \\
(II) $L\in\mathcal X$ whenever
$L\le G$ contains some $H\in\mathcal X$ as a subgroup of finite index; \\
(III)
$H\in\mathcal X$ whenever $H\le\St_G(1)$ is finitely generated and
$\{h|_x:h\in H\}\in\mathcal X$ for every $x\in X$. 

A finitely generated group $G$
has the \emph{subgroup induction property} if every inductive collection contains
all finitely generated subgroups of $G$; the weaker version used in \cite{GLN21} is
called the weak subgroup induction property in \cite{FL25,FGLN25}. Moreover, $G$
acts \emph{tree-primitively} if the only $G$-invariant partitions of $\Lv_n$ are
the partitions $\{vX^{n-m}:v\in X^m\}$, $0\le m\le n$ \cite[Definition~4.1]{FGLN25}.

For a non-empty finite set $V$ of pairwise incomparable vertices, a \emph{diagonal
block subgroup} with supporting set $V$ is a subgroup
$D\le\bigcap_{v\in V}\St_G(v)$ fixing every vertex outside $\bigcup_{v\in V}\T_v$
such that, for each $v\in V$, the map $g\mapsto g|_v$ is injective on $D$ with image
of finite index in $\{g|_v:g\in\St_G(v)\}$. A \emph{block subgroup} is a product of
diagonal block subgroups with disjoint supporting sets whose union consists of
pairwise incomparable vertices \cite[Definitions~5.1 and~5.4]{FGLN25}. For instance,
if $[G:K]<\infty$, the subgroups $P$ and $Q$ of Theorem~\ref{thm:A} are diagonal
block subgroups supported on $\{0,(d-1)(d-1)\}$ and $\{0,d-1\}$ (compare
\cite[Example~5.3]{FGLN25}). By \cite[Theorem~1.1]{FGLN25}, every finitely generated
subgroup of a finitely generated self-replicating branch group with the subgroup
induction property acting tree-primitively contains a block subgroup of finite
index; a version for the weak property was announced in \cite[Theorem~2.9]{GLN21}
(see \cite[Remark~5.9]{FGLN25}).

\begin{corollary}\label{cor:hered}
Let $G$ satisfy the hypotheses of Theorem~\ref{thm:A}, and suppose that it is
self-replicating, has the subgroup induction property and acts tree-primitively.
Then no infinite finitely generated subgroup of $G$ is double coset separable. This
applies to the first Grigorchuk group and to every torsion GGS-group, whose finitely
generated subgroups are moreover subgroup separable.
\end{corollary}

\begin{proof}
Let $L\le G$ be infinite and finitely generated. By \cite[Theorem~1.1]{FGLN25}, $L$
has a subgroup of finite index $D_1\cdots D_r$, with $r\ge1$ and each $D_i$ a
diagonal block subgroup. Since $G$ is self-replicating, $D_1$ is isomorphic to a
subgroup of finite index of $G$; by Lemma~\ref{lem:hered} and Theorem~\ref{thm:A},
neither $D_1$ nor $L$ is double coset separable. The first Grigorchuk group and the
torsion GGS-groups satisfy the hypotheses of Theorem~\ref{thm:A} (a torsion
GGS-group has non-constant defining vector); they are self-replicating, have the
subgroup induction property \cite{GW03,FL25}, act tree-primitively
\cite[Corollaries~4.6 and~4.7]{FGLN25}, and are subgroup separable
\cite[Theorem~C]{FL25}, a property inherited by subgroups.
\end{proof}
\section{The Basilica group and the \texorpdfstring{$p$}{p}-Basilica groups}
\label{sec:basilica}

In this section we apply Theorem~\ref{thm:C} to the Basilica group and to its
generalisations to trees of prime degree. For a prime $p$, the
\emph{$p$-Basilica group} is $\B_p=\langle a,b\rangle\le\Aut(T_p)$, where
\[
 a=(1,\ldots,1,b),\qquad
 b=(1,\ldots,1,a)\sigma,\qquad
 \sigma=(0\ 1\ \cdots\ p-1).
\]
The group $\B_2$ is the Basilica group. These groups are self-similar,
level-transitive and weakly regular branch over their derived subgroups
\cite[Lemma~3.1 and Theorem~3.5]{DFNT22}. Each generator has exactly one non-trivial
section at every level, so that they have bounded activity. On the other hand, they
are not branch groups and do not have the congruence subgroup property
\cite[Theorems~A and~B]{DFNT22}. The remaining hypothesis of Theorem~\ref{thm:C} is
the finite generation of the derived subgroup. Francoeur established it for $\B_2'$
\cite[Proposition~4.7]{Fra20}; we now prove it for every prime $p$.

\begin{lemma}\label{lem:pbasilica-derived}
For every prime $p$, the derived subgroup $\B_p'$ is generated by the following
$2p-1$ commutators:
\[
 [a,b^p],\qquad [a,b^r],\qquad [a,b^{r-p}]
 \quad(1\le r<p).
\]
\end{lemma}

\begin{proof}
Let us put $z_t=[a,b^t]$ for $t\in\mathbb Z$, and let $L$ be the subgroup generated
by the displayed elements. Since $\B_p$ is generated by $a$ and $b$, its derived
subgroup is the normal closure of $z_1$. The elements $a$ and $a^b$ are supported in
distinct first-level subtrees. Thus $c=[b,a]=(a^b)^{-1}a$ commutes with $a$, and
$b^{a^m}=bc^m$ for every $m\in\mathbb Z$.

Since $b^p=(a,\ldots,a)$, the automorphism $b^{pm+r}$ sends $p-1$ to $r-1$
and has section $a^{m+1}$ at $p-1$, for $m\in\mathbb Z$ and $1\le r<p$.
Thus $a^{b^{pm+r}}=\lambda_{r-1}(b^{a^{m+1}})$, whereas
$a^{b^{pm}}=\lambda_{p-1}(b^{a^m})$. Using $b^{a^m}=bc^m$, we obtain
\[
 z_{pm}=\lambda_{p-1}(c^m),\qquad
 z_{pm+r}=\lambda_{p-1}(b^{-1})\,\lambda_{r-1}(bc^{m+1}).
\]
The two factors of $z_{pm+r}$ commute, and
$z_{r-p}^{-1}z_r=\lambda_{r-1}(c)$. Consequently,
\[
 z_{pm}=z_p^m,\qquad
 z_{pm+r}=z_r\bigl(z_{r-p}^{-1}z_r\bigr)^m,
\]
so that $z_t\in L$ for every integer $t$.

Let us show that $L$ is normal in $\B_p$. For $\varepsilon\in\{-1,1\}$, the
commutator identity $z_t^{b^\varepsilon}=z_\varepsilon^{-1}z_{t+\varepsilon}$ shows
that $b$ normalises $L$. If $p\nmid t$, then $z_t$ commutes with $a$, since $a$ and
$a^{b^t}$ have disjoint supports. If $p\mid t$, then $z_t$ is supported below the
vertex $p-1$, where $a$ and $z_1^{-1}$ have the same section $b$; hence
$z_t^{a^\varepsilon}=z_t^{z_1^{-\varepsilon}}\in L$. Thus $a$ normalises $L$ as
well, and $L\trianglelefteq\B_p$. Since $L$ contains $z_1=[a,b]$, whose normal
closure is $\B_p'$, and is contained in $\B_p'$, we conclude that $L=\B_p'$.
\end{proof}

For $p=2$, Lemma~\ref{lem:pbasilica-derived} gives exactly the generating set
$\{[a,b],[a,b^{-1}],[a,b^2]\}$ of \cite[Proposition~4.7]{Fra20}, which is minimal
since $\B_2'/\B_2''\cong\mathbb Z^3$ \cite[Proposition~4.14]{Fra20}. We can now
prove the announced result.

\begin{corollary}\label{cor:basilica}
For every prime $p$, the group $\B_p$ contains closed finitely generated subgroups
$P,Q\cong\B_p'$ with $P\cap Q=1$ and such that $PQ$ is not closed. In particular,
neither the Basilica group nor any $p$-Basilica group is double coset separable.
\end{corollary}

\begin{proof}
The subgroup $K=\B_p'$ is finitely generated by Lemma~\ref{lem:pbasilica-derived},
and it is closed since $\B_p/\B_p'\cong\mathbb Z^2$ is residually finite
\cite[Theorem~A]{DFNT22}. Hence Theorem~\ref{thm:C}, including its final assertion,
applies.
\end{proof}

\section{Polynomial activity and the long-range group}\label{sec:linear}
In the previous sections we used the bounded activity. Now we want to explore how much such a hypothesis is restrictive. 
The activity calculation also applies to polynomial activity, provided that the
chosen element has a sufficiently large activity function. We first state the
criterion, then verify its hypotheses for the long-range group.

\subsection{Polynomial activity}

\begin{theorem}\label{thm:Cpoly}
Let $G\le\Aut(T_d)$ be level-transitive, and let $m\ge0$ be a real number such that
$A_n(g)=O(n^m)$ for every $g\in G$. Suppose that $1\ne K\le G$ is finitely
generated and that $\lambda_x(K)\le K$ for every $x\in X$, and define $P$ and $Q$ as
in Theorem~\ref{thm:A}. Then $P\cong Q\cong K$ and $P\cap Q=1$. Moreover, every
$k\in K'$ with $\liminf_{n\to\infty}A_n(k)/n^m>0$ satisfies
\[
 \lambda^2(k)\in\cl{PQ}\setminus PQ .
\]
Hence, if such an element exists, $G$ is not double coset separable. If, in
addition, $G$ is self-similar and $K$ is closed in $G$, then $P$ and $Q$ are closed
in $G$.
\end{theorem}

\begin{proof}
For $m=0$ the assumption on $k$ says that $k$ is non-finitary, and we recover
Theorem~\ref{thm:C}. The proofs of the algebraic assertions and of membership in
$\cl{PQ}$ are unchanged, as is the proof of closedness of the factors. It remains
to show that $\lambda^2(k)\notin PQ$. If $\lambda^2(k)\in PQ$,
then, as in the proof of Theorem~\ref{thm:A}, there is $x\in K$ with
$x=\lambda(k)\rho(x)$, and the recursion in the proof of Lemma~\ref{lem:activity}
gives $A_n(x)=1+\sum_{j<n}A_j(k)$. By assumption there are $c>0$ and $n_0$ such that
$A_j(k)\ge cj^m$ for $j\ge n_0$. If $n\ge\max\{2n_0,2\}$, then the
$\lfloor n/2\rfloor\ge n/4$ indices $j$ with $n/2\le j<n$ satisfy $j\ge n_0$, hence
\[
 A_n(x)\ge\frac n4\cdot c\Bigl(\frac n2\Bigr)^m=\frac{c}{2^{m+2}}\,n^{m+1},
\]
which contradicts $A_n(x)=O(n^m)$, since $x\in G$.
\end{proof}

\subsection{The long-range group}

In this subsection $\sigma=(0\ 1)$ and $G=\langle a,b\rangle\le\Aut(T_2)$, where
\[
 a=(1,a)\sigma,\qquad b=(a,b).
\]
The element $a$ is the binary odometer: it adds one to a binary word with the
first letter as the least significant digit. Thus it acts transitively on every
level. The displayed recursions also show that $G$ is self-similar. For $t\in\mathbb Z$ we put $b_t=b^{a^t}$, and we use the commutators
\[
 z_t=[b,a^t]=b^{-1}b_t,\qquad y_t=[b^{-1},a^t]=b\,b_t^{-1}.
\]
For words in $a^{\pm1},b^{\pm1}$, we write $e_a,e_b$ for their exponent sums.
At this point they are functions of words; Lemma~\ref{lem:lin-ab} will show that
they depend only on the represented element of $G$. We put $G^{ab}=G/G'$ and
write $\bar g$ for the image of $g$ in this quotient.

\begin{remark}\label{rem:long-range}
Identify the boundary $X^\omega$, the set of infinite binary words, with the ring
$\mathbb Z_2$ of $2$-adic integers by
$x_0x_1x_2\cdots\mapsto1+\sum_{i\ge0}x_i2^i$. Then $a$ acts as $u\mapsto u+1$.
The point $1^\omega$ corresponds to $0$ and is fixed by $b$. Every other point
has the form $1^j0\eta$; here $b$ fixes the prefix $1^j0$ and applies $a$ to
$\eta$. Hence $b$ acts as $u\mapsto u+2^{\nu_2(u)+1}$ for $u\ne0$, where
$\nu_2(u)$ is the exponent of $2$ dividing $u$. Thus $b$ adds the smallest power
of $2$ which does not divide $u$. The subset $\mathbb Z\subset\mathbb Z_2$
is the orbit of $0$, which corresponds to the point $1^\omega$, and it is invariant and
dense, so that $G$ acts faithfully on it. Up to interchanging the two generators, this
is the description of the long-range group given in \cite[Section~1]{AAV13}; its Schreier
graph on $\mathbb Z$ was studied in \cite{BH05, BCDN12}. Hence $G$ is isomorphic to the
long-range group. Moreover, if $s=(s,s)\sigma$ is the automorphism exchanging the
letters $0$ and $1$, then $a^s=(a^s,1)\sigma$ and $b^s=(b^s,a^s)$, so that
$\alpha=(a^s)^{-1}$ and $\beta=(b^s)^{-1}$ satisfy $\alpha=(1,\alpha)\sigma$ and
$\beta=(\beta,\alpha)$. Thus $G$ is conjugate in $\Aut(T_2)$ to the group
$\langle\alpha,\beta\rangle$, which is another common presentation of the long-range
group.
\end{remark}

\begin{lemma}\label{lem:lin-activity}
We have $A_n(a)=1$ and $A_n(b)=n+1$ for every $n\ge0$. Consequently
$A_n(g)\le(n+1)|g|$ for every $g\in G$, where $|g|$ is the word length with respect
to $\{a,b\}$.
\end{lemma}

\begin{proof}
The only non-trivial section of $a$ on $\Lv_n$ is $a$, at $1^n$. The non-trivial
sections of $b$ on $\Lv_n$ are $b$ at $1^n$ and $a$ at $1^j0\,1^{n-j-1}$ for
$0\le j<n$. The last assertion follows from the subadditivity of $A_n$.
\end{proof}

\begin{lemma}\label{lem:lin-ab}
The abelianisation of $G$ is free abelian of rank $2$, with basis $\bar a,\bar b$.
\end{lemma}

\begin{proof}
Let $F$ be the free group on $a,b$, let $R\trianglelefteq F$ be the kernel of
$F\to G$, and let $\Lambda\le\mathbb Z^2$ be the image of $R$ under
$w\mapsto(e_a(w),e_b(w))$. We have to show that $\Lambda=0$.

For a word $w$ and $x\in\{0,1\}$, read $w$ from left to right along the walk on
$\{0,1\}$ which starts at $x$, is switched by each letter $a^{\pm1}$ and is fixed by
each letter $b^{\pm1}$. Let $w|_x$ be the word obtained by recording $a$ for each $a$
read at position $1$, $a^{-1}$ for each $a^{-1}$ read at position $0$, $a^{\pm1}$ for
each $b^{\pm1}$ read at position $0$, and $b^{\pm1}$ for each $b^{\pm1}$ read at
position $1$, recording the empty word in the remaining cases. By
\eqref{eq:sections}, $w|_x$ represents the section at $x$ of the
element represented by $w$; in particular, $w\in R$ implies $w|_x\in R$.

Let $w\in R$ and $(A,B)=(e_a(w),e_b(w))$. The root permutation of $w$ is $\sigma^A$,
so $A$ is even and both walks are closed. Comparing the numbers of transitions
$0\to1$ and $1\to0$ of a closed walk, one finds that the number of letters $a$ read
at position $1$, minus the number of letters $a^{-1}$ read at position $0$, equals
$A/2$. Since the two walks occupy complementary positions at every step, putting
$\beta=e_b(w|_1)$ we obtain
\[
 (e_a,e_b)(w|_0)=\Bigl(\tfrac A2+\beta,\ B-\beta\Bigr),\qquad
 (e_a,e_b)(w|_1)=\Bigl(\tfrac A2+B-\beta,\ \beta\Bigr).
\]
In particular $(B,0)=(e_a,e_b)(w|_0)+(e_a,e_b)(w|_1)-(A,B)\in\Lambda$.

We show by induction on $j\ge1$ that $2^j$ divides both coordinates of every element
of $\Lambda$. For $j=1$, the first coordinate of every element of $\Lambda$ is even,
and so is $B$, being the first coordinate of $(B,0)\in\Lambda$. Assume the claim for
$j$ and let $(A,B)\in\Lambda$, coming from $w\in R$. Applying the claim to
$(e_a,e_b)(w|_1)$ and $(e_a,e_b)(w|_0)$ gives $2^j\mid\beta$ and $2^j\mid A/2+\beta$,
hence $2^{j+1}\mid A$. Applying this to $(B,0)\in\Lambda$ gives $2^{j+1}\mid B$.
Therefore $\Lambda=0$.
\end{proof}

Thus $e_a,e_b$ define homomorphisms $G\to\mathbb Z$, which we use below.

\begin{corollary}\label{cor:lin-structure}
\begin{enumerate}
\item[(a)] The subgroups $\psi(\St_G(1))$ and $\psi(\rist_G(1))$ have infinite index
in $G\times G$. Hence $G$ is not regular branch over any subgroup of finite index and
is not a branch group. Moreover, $G$ does not have the congruence subgroup property.
\item[(b)] The assignment $a\mapsto a^2$, $b\mapsto b$ extends to an injective
endomorphism $\varphi$ of $G$ with $\psi(\varphi(g))=(a^{e_a(g)+e_b(g)},g)$.
\item[(c)] For $g\in G'$ we have $\rho(g)=\varphi(g)$ and
$\lambda(g)=a^{-1}\rho(g^{a^{-1}})a$. In particular $\lambda(G')\le G'$ and
$\rho(G')\le G'$, so that $G$ is weakly regular branch over $G'$.
\end{enumerate}
\end{corollary}

\begin{proof}
(a) By Schreier's lemma, $\St_G(1)=\langle b,b^a,a^2\rangle$, and
$\psi(b)=(a,b)$, $\psi(b^a)=(b^a,a)$, $\psi(a^2)=(a,a)$. By
Lemma~\ref{lem:lin-ab}, the image of $\psi(\St_G(1))$ in
$G^{ab}\times G^{ab}\cong\mathbb Z^4$ is spanned by $(1,0,0,1)$, $(0,1,1,0)$ and
$(1,0,1,0)$, so it has rank $3$ and therefore infinite index in $\mathbb Z^4$.
This excludes regular branching over a finite-index subgroup $K$, since it would
give $K\times K\le\psi(K\cap\St_G(1))$ of finite index in $G\times G$.
An element $(x,1)$ or $(1,x)$ of $\psi(\St_G(1))$ has
$\bar x\in\mathbb Z(-1,1)$. Consequently, the image of $\psi(\rist_G(1))$ has
rank at most $2$, strictly less than the rank-$3$ image of $\psi(\St_G(1))$.
It follows that $[\St_G(1):\rist_G(1)]=\infty$, so $G$ is not branch.
Finally, all congruence quotients are $2$-groups, whereas the abelianisation
provides a quotient $C_3$. Its kernel cannot contain any level stabiliser.

(b) Let $F=F(a,b)$, and let $\tilde\varphi\colon F\to F$ be defined by
$a\mapsto a^2$, $b\mapsto b$. Since $\psi(a^{\pm2})=(a^{\pm1},a^{\pm1})$ and $\psi(b^{\pm1})=(a^{\pm1},b^{\pm1})$, the
element represented by $\tilde\varphi(w)$ has first-level sections
$a^{e_a(w)+e_b(w)}$ and the element represented by $w$. If $w$ represents the
identity in $G$, then $e_a(w)=e_b(w)=0$ by Lemma~\ref{lem:lin-ab}, so
$\tilde\varphi(w)$ also represents the identity. Thus
$\tilde\varphi$ induces an endomorphism $\varphi$ of $G$, which is injective since
$g$ is a section of $\varphi(g)$.

(c) If $g\in G'$, then $\psi(\varphi(g))=(1,g)$, that is, $\varphi(g)=\rho(g)$, and
$\varphi(G')\le G'$. The formula for $\lambda$ follows from $a^{-1}(1,h)a=(h^a,1)$,
and $G'$ is normal in $G$. Also, $G'\ne1$: the section
$[b,a]|_0=a^{-2}ba$ has non-trivial root permutation. Thus the branching subgroup
is non-trivial, as required by Definition~\ref{def:branch}.
\end{proof}

\begin{lemma}\label{lem:lin-relation}
For $m\in\mathbb Z$ we have $b_{2m}=(a,b_m)$ and $b_{2m+1}=(b_{m+1},a)$. Moreover,
\[
 b_t\,b_{t+1}=b_{t-1}\,b_{t+2}\qquad(t\in\mathbb Z).
\]
Consequently, if $c=ba^{-1}$ and $c_t=c^{a^t}$, then $c_{t+1}=c_{t-1}^{-1}c_t^2$ for
every $t$, and the normal closure $N$ of $c$ in $G$ is generated by $c_0$ and $c_1$.
\end{lemma}

\begin{proof}
The section formula gives $(x,y)^a=(y^a,x)$. Since
$a^{2m}=(a^m,a^m)$ for every $m\in\mathbb Z$, conjugating $b=(a,b)$ first by
$a^{2m}$ and then by $a$ gives the first two formulas, including negative $m$. Since $b_2=b_1^a$ and
$b_1=b^a$, we get $b_1b_2=(b_1a,ab_1)=(ab_2,ba)=b_0b_3$, and conjugating by
$a^{t-1}$ gives the relation. Since $b_t=c_ta$ and $ac_sa^{-1}=c_{s-1}$, the relation
becomes $c_t^2=c_{t-1}c_{t+1}$. Finally, the subgroup generated by the elements $c_t$
is normalised by $a$ and contains $c$, and $G=\langle a,c\rangle$; hence it equals
$N$, and by the recursion it is generated by $c_0$ and $c_1$.
\end{proof}

\begin{proposition}\label{prop:lin-fg}
The derived subgroup $G'$ is generated by
\[
 [b,a],\qquad [b^{-1},a^2],\qquad [b^{-1},a]^{-1}[b^{-1},a^3],\qquad
 [b^{-1},a^{-1}]^{-1}[b^{-1},a^3].
\]
\end{proposition}

\begin{proof}
Recall that $c_t=a^{-t}ba^{t-1}$, so that $c_t^a=c_{t+1}$, and that
$y_t=[b^{-1},a^t]=c_0c_t^{-1}$. For $x\in\St_G(1)$ with $\psi(x)=(x_0,x_1)$ we have
\begin{equation}\label{eq:lin-conj}
 \psi(x^a)=(x_1^a,x_0).
\end{equation}

Step 1: $G'$ as a normal closure inside $N$. In fact, in $G/N$ the elements $a$ and
$b$ have the same image, so $G/N$ is cyclic; since $e_a+e_b$ vanishes on $N$ and
sends $a$ to $1$, we get $G/N\cong\mathbb Z$ and $N=\ker(e_a+e_b)$. By
Lemma~\ref{lem:lin-ab}, $N$ is the inverse image of $\mathbb Z\cdot(-1,1)$ in
$G^{ab}$; hence $G'\le N$ and $N/G'$ is infinite cyclic, generated by the image of
any $c_t$. We use the following elementary fact: \emph{if $H=\langle
g_1,\ldots,g_r,s\rangle$ with $g_i\in L\trianglelefteq H$ and $H/L$ infinite cyclic
generated by $sL$, then $L=\langle g_i^{s^k}:k\in\mathbb Z\rangle$.} Indeed, the
latter subgroup $L_0$ is normalised by $s^{\pm1}$ and by the $g_i$, so
$H=L_0\langle s\rangle$ and $L=L_0(L\cap\langle s\rangle)=L_0$, since
$L\cap\langle s\rangle=1$. Put $s_1=c_2$ and $s_2=c_1^{-1}$. Since
$N=\langle c_0,c_1\rangle$ (Lemma~\ref{lem:nonfin}), we have
$N=\langle c_0c_2^{-1},c_1c_2^{-1},s_1\rangle=\langle c_0c_1^{-1},s_2\rangle$, and
therefore
\begin{equation}\label{eq:lin-Gprime}
 G'=\bigl\langle (c_0c_2^{-1})^{s_1^k},(c_1c_2^{-1})^{s_1^k}:k\in\mathbb Z\bigr\rangle
   =\bigl\langle (c_0c_1^{-1})^{s_2^k}:k\in\mathbb Z\bigr\rangle .
\end{equation}

Step 2: determination of the image of $G'$ under $\psi$. Since $[G:\St_G(1)]=2$, we have
$G'\le\St_G(1)$, and a direct computation gives
$\theta:=\psi(z_1)=(s_1,s_2)$. Let $q\colon G\times G\to G^{ab}\times G^{ab}$ be the
product of the abelianisation maps; its kernel $G'\times G'$ is contained in
$\psi(G')$ by Corollary~\ref{cor:lin-structure}(c). Since $G'$ is the normal
closure of $z_1$, and every element of $G$ is of the form $h$ or $ah$ with
$h\in\St_G(1)$, \eqref{eq:lin-conj} shows that $\psi(G')$ is generated by the
elements
\[
 \psi(z_1^h)=(s_1^{h_0},s_2^{h_1}),\qquad
 \psi(z_1^{ah})=(s_2^{ah_0},s_1^{h_1}),\qquad \psi(h)=(h_0,h_1).
\]
As $\bar s_2=-\bar s_1$, their images under $q$ are $\pm q(\theta)$. Hence
$q(\psi(G'))=\langle q(\theta)\rangle$ and, since $\ker q\le\psi(G')$,
\[
 \psi(G')=(G'\times G')\langle\theta\rangle .
\]
Since $(u,1)^{\theta^k}=(u^{s_1^k},1)$ and $(1,u)^{\theta^k}=(1,u^{s_2^k})$,
\eqref{eq:lin-Gprime} shows that $\psi(G')$ is generated by $\theta$,
$(c_0c_2^{-1},1)$, $(c_1c_2^{-1},1)$ and $(1,c_0c_1^{-1})$, that is,
\[
 G'=\bigl\langle z_1,\ \lambda(y_2),\ \lambda(y_1^a),\ \rho(y_1)\bigr\rangle,
\]
where we used $c_0c_1^{-1}=y_1$, $c_0c_2^{-1}=y_2$ and
$c_1c_2^{-1}=(c_0c_1^{-1})^a=y_1^a$.

\emph{Step 3: explicit generators.} From $b_{2t}=(a,b_t)$ we get
$\rho(y_t)=y_{2t}$, and the commutator identity $[x,yz]=[x,z][x,y]^z$ gives
$y_t^{a^{\pm1}}=y_{\pm1}^{-1}y_{t\pm1}$. Using also \eqref{eq:lin-conj}, we obtain
\[
 \begin{aligned}
 \rho(y_1)&=y_2,\\
 \lambda(y_1^a)&=y_2^a=y_1^{-1}y_3
   &&\text{since }\psi(y_2)=(1,y_1),\\
 \lambda(y_2)&=\rho\bigl(y_2^{a^{-1}}\bigr)^a
   =\rho\bigl(y_{-1}^{-1}y_1\bigr)^a=(y_{-2}^{-1}y_2)^a=y_{-1}^{-1}y_3 .
 \end{aligned}
\]
Since $z_1=[b,a]$ and $y_t=[b^{-1},a^t]$, these are the four generators in the
statement.
\end{proof}

\begin{lemma}\label{lem:lin-witness}
The element $k=[z_1,\lambda_{00}(z_1)]$ belongs to $G''$, and $A_n(k)\ge n-2$ for
every $n\ge2$.
\end{lemma}

\begin{proof}
By Corollary~\ref{cor:lin-structure}(c), $\lambda_{00}(z_1)=\lambda^2(z_1)\in G'$, so
$k\in G''$. The section of $z_1$ at $0$ is $s_1=a^{-2}ba$, which acts as $\sigma$ on
the first level; hence $z_1$ maps $00$ to $01$, and
$k=\bigl(\lambda_{00}(z_1)^{-1}\bigr)^{z_1}\lambda_{00}(z_1)$ is the product of an
element supported on $\T_{01}$ and of $\lambda_{00}(z_1)$. Therefore
$A_n(k)\ge A_{n-2}(z_1)$. The section of $z_1$ at $1$ is
$s_2=b^{-1}a=(a^{-1},b^{-1}a)\sigma$, whose non-trivial sections on $\Lv_m$ are
$b^{-1}a$ at $1^m$ and $a^{-1}$ at $1^j0^{m-j}$ for $0\le j<m$. Hence
$A_m(z_1)\ge A_{m-1}(s_2)=m$ for $m\ge1$.
\end{proof}

\begin{proof}[Proof of Theorem~\ref{thm:D}]
Remark~\ref{rem:long-range} identifies $G$ with the long-range group.
Lemma~\ref{lem:lin-activity} gives linear activity, and
Corollary~\ref{cor:lin-structure} shows that $G$ is weakly regular branch over
$G'$ but is not branch. By Proposition~\ref{prop:lin-fg}, $G'$ is finitely
generated. It is closed, since $G/G'\cong\mathbb Z^2$ is residually finite.

Apply Theorem~\ref{thm:Cpoly} with $m=1$ and $K=G'$. The element
$k=[z_1,\lambda_{00}(z_1)]\in G''$ satisfies $A_n(k)\ge n-2$ by
Lemma~\ref{lem:lin-witness}, so
\[
 \lambda^2(k)\in\cl{PQ}\setminus PQ.
\]
The same theorem gives $P\cong Q\cong G'$, $P\cap Q=1$ and the closedness of
both factors.
\end{proof}

\begin{remark}\label{rem:linear-witness}
Non-finitarity alone does not rule out solutions of the recursive equation in a
group of linear activity. Indeed, $\psi(b)=(a,b)$ means that $b=\lambda(a)\rho(b)$: the equation
of Lemma~\ref{lem:activity}, with the bounded non-finitary element $a$ in place of
$k$, has the solution $b\in G$. Similarly
$[a,b^{-1}]=[b,a^{-1}]\,\rho([a,b^{-1}])$, so that a twisted equation $x=c\rho(x)$
with coefficient in $G'$ has a solution in $G'$. Theorem~\ref{thm:B} does not
apply, since $G$ is not branch. The closure of $PQ$ remains undetermined.
\end{remark}

\section{The Hanoi towers example}\label{sec:hanoi}

The approximation in the proof of Theorem~\ref{thm:A} uses an element of $K'$. In
this section we show, using the three-peg Hanoi towers group, that a non-finitary
element of $K$ need not give a point of $\cl{PQ}$. Let
$H=H^{(3)}=\langle a_1,a_2,a_3\rangle\le\Aut(T_3)$, where
\[
 a_1=(a_1,1,1)(1\ 2),\qquad
 a_2=(1,a_2,1)(0\ 2),\qquad
 a_3=(1,1,a_3)(0\ 1).
\]
This group is self-similar, level-transitive and regular branch over $H'$
\cite{GS06,BSZ12}. Its generators satisfy $A_n(a_i)=1$, so that
Theorem~\ref{thm:A} applies with $K=H'$. Notice that $H'$ does not fix the first
level; this is the reason why Definition~\ref{def:branch} involves
$K\cap\St_G(1)$.

\begin{proposition}\label{prop:hanoi}
Let $K=H'$ and define $P,Q$ as in Theorem~\ref{thm:A}. If
$k\in H'\setminus H''(H')^2$, where $(H')^2$ is the subgroup generated by the
squares of the elements of $H'$, then $\lambda^2(k)\notin\cl{PQ}$. This applies to
the non-finitary element $k=[a_1,a_2]$.
\end{proposition}

\begin{proof}
By \cite[Proposition~3.17]{BSZ12}, $H'/H''\cong\mathbb Z^3\times C_3$, where $C_3$
is cyclic of order $3$, and the images of $[a_1,a_2],[a_2,a_3],[a_3,a_1]$ generate
the free abelian factor. Thus $V:=H'/H''(H')^2\cong(\mathbb Z/2\mathbb Z)^3$. The
calculation in \cite[proof of Proposition~3.18]{BSZ12} shows that $\lambda_0$ induces
$\operatorname{diag}(1,-1,1)$ on the free abelian factor, up to terms in $C_3$; hence
it induces the identity on $V$. Let $s$ be the automorphism acting at every vertex
as the cycle $(0\ 1\ 2)$. It normalises $H$ and permutes its generators cyclically.
The identity $\lambda_{i^s}(x^s)=\lambda_i(x)^s$, together with the invariance of
$H''(H')^2$ under conjugation by $s$, shows that the maps induced by the three
rigid embeddings on $V$ are conjugate. Thus $\rho=\lambda_2$ also induces the
identity. In particular,
$\overline{\lambda(x)}=\overline{\rho(x)}=\bar x$ for every $x\in K$.

Let $\Phi\colon W'=\psi^{-1}(K^3)\to V^3$ be the section map followed by the
reduction modulo $H''(H')^2$ in each coordinate. Writing $\bar x$ for the image of
$x\in K$ in $V$, we obtain
\[
 \Phi(PQ)\subseteq\{(u,0,u):u\in V\},\qquad
 \Phi(\lambda^2(k))=(\bar k,0,0),
\]
and the latter vector does not belong to the former set when $\bar k\ne0$. The kernel
of $\Phi$ has finite index in $W'$, hence in $H$, and its normal core is a normal
subgroup of finite index $N\trianglelefteq H$ contained in $\ker\Phi$. The same
calculation therefore separates $\lambda^2(k)$ from $PQ$ in $H/N$, that is,
$\lambda^2(k)\notin\cl{PQ}$.

Finally, a direct computation gives $[a_1,a_2]=(a_1a_2,a_1,a_2)(0\ 1\ 2)$. Its
section at $1$ is $a_1$, and $a_1|_0=a_1\ne1$; thus $[a_1,a_2]$ is non-finitary, and
its image in $V$ is non-zero by the description of the free abelian factor.
\end{proof}

\begin{remark}
Since $\lambda^2(k)\in W$, Proposition~\ref{prop:hanoi} and Theorem~\ref{thm:B} show
in particular that $H^{(3)}$ does not have trivial branch kernel (compare
\cite{BSZ12}).
\end{remark}

\section{Further consequences}\label{sec:twisted}

\subsection{Twisted conjugacy classes}
The construction can be interpreted in terms of twisted conjugacy classes. We record
the general correspondence for an arbitrary group $H$. For an endomorphism
$\alpha\in\End(H)$, the \emph{$\alpha$-twisted conjugacy class} of $g\in H$ is
$[g]_\alpha=\{xg\alpha(x)^{-1}:x\in H\}$. In particular, putting
$D_\alpha(x)=x\alpha(x)^{-1}$, we have $[1]_\alpha=D_\alpha(H)$. Let
$\Gamma_\alpha=\{(x,\alpha(x)):x\in H\}$ and $\Delta_H=\{(y,y):y\in H\}$. The
following proposition explains why these classes appear in our construction.

\begin{proposition}\label{prop:twisted}
For every $g\in H$,
\[
 \Gamma_\alpha(g,1)\Delta_H
   =\{(u,v)\in H\times H:uv^{-1}\in[g]_\alpha\}.
\]
Its closure in $H\times H$ is $\{(u,v):uv^{-1}\in\cl{[g]_\alpha}\}$, where the
closure on the right is taken in $H$. Thus the double coset
$\Gamma_\alpha(g,1)\Delta_H$ is closed if and only if $[g]_\alpha$ is closed.
Moreover, $\Gamma_\alpha\cap\Delta_H\cong\Fix(\alpha)$, where
$\Fix(\alpha)=\{x\in H:\alpha(x)=x\}$.
\end{proposition}

\begin{proof}
We have
\[
 (x,\alpha(x))(g,1)(y,y)=(xgy,\alpha(x)y),\qquad
 (xgy)(\alpha(x)y)^{-1}=xg\alpha(x)^{-1};
\]
conversely, if $uv^{-1}=xg\alpha(x)^{-1}$, it is enough to take
$y=\alpha(x)^{-1}v$.

The homeomorphism $(u,v)\mapsto(uv^{-1},v)$, whose inverse is
$(z,v)\mapsto(zv,v)$, sends the displayed double coset onto $[g]_\alpha\times H$.
It therefore sends its closure onto $\cl{[g]_\alpha}\times H$, proving the closure
formula and the equivalence of closedness. Finally, $(x,\alpha(x))$ belongs to the diagonal exactly
when $\alpha(x)=x$.
\end{proof}

\begin{corollary}\label{cor:tcs}
Let $H$ be finitely generated and residually finite. If $H\times H$ embeds in a
double coset separable group, then $[g]_\alpha$ is closed for every $g\in H$ and
every $\alpha\in\End(H)$. In particular, $H$ is conjugacy separable, that is, every
conjugacy class of $H$ is closed.
\end{corollary}

\begin{proof}
By Lemma~\ref{lem:hered}, $H\times H$ is double coset separable. The groups
$\Gamma_\alpha$ and $\Delta_H$ are finitely generated, so that the double coset
$\Gamma_\alpha(g,1)\Delta_H$ is closed and Proposition~\ref{prop:twisted} applies.
For the last assertion, take $\alpha$ to be the identity.
\end{proof}

\subsection{Rational subsets}
For the first Grigorchuk group, our construction also gives a non-closed rational
subset. Recall that a subset of a finitely generated group is \emph{rational} if it
is the image under evaluation of a language accepted by a finite automaton; equivalently,
rational subsets are obtained from finite subsets by finite unions, products and
taking the submonoid generated by a subset. Finitely generated subgroups are
rational, so that $R=PQ\lambda^2((ab)^2)^{-1}$ is rational and
$1\in\cl R\setminus R$.

In a torsion group, every submonoid is a subgroup; if it is rational, then it is
finitely generated by the Anisimov--Seifert theorem \cite{AS75}. Induction on
rational expressions then shows that every rational subset of a torsion group is a
finite union of translates of finite products of finitely generated subgroups. Thus,
in a residually finite torsion group, all rational subsets are closed if and only if
the group is product separable.

\subsection{Deciding membership in the Grigorchuk example}
Since the product $PQ$ of Section~\ref{sec:examples} is not closed, finite
quotients cannot separate every element outside $PQ$ from $PQ$. Nevertheless,
membership in this product is decidable.

\begin{proposition}\label{prop:membership}
Membership in the product $PQ$ constructed in Section~\ref{sec:examples} for the
first Grigorchuk group is decidable.
\end{proposition}

\begin{proof}
Recall that the generators of $\Gr$ are finite-state. By \eqref{eq:sections} we can
compute automata for products of
finite-state automorphisms, and equality of two finite-state automorphisms is
decidable by comparing root permutations at every reachable pair of states. The
quotient $\Gr/K$ is effectively computable as well: computing, in the quotients
$\Gr/\St_\Gr(n)$, the normal closure of the image of $(ab)^2$ until its index is
$16$, we find $n_0$ with $\St_\Gr(n_0)\le K$, and the congruence subgroup property
guarantees that this procedure terminates. In particular, for $h\in\Gr$ the coset
$hK$ only depends on the action of $h$ on $\Lv_{n_0}$.

Let $g\in\Gr$ be given by a word. We first
check that $g$ fixes the first level and that its sections $u=g|_0$ and $v=g|_1$
belong to $K$; otherwise $g\notin W\supseteq PQ$. For $g\in W$, the defining
equations of $P$ and $Q$ give
\[
 g\in PQ\quad\Longleftrightarrow\quad
 x=z\rho(x)\text{ has a solution }x\in K,\qquad z=uv^{-1}.
\]
The equation $x=z\rho(x)$ determines a unique automorphism $x\in\Aut(T_2)$, since
its action on each finite level is determined recursively. We construct a finite
automaton for $x$ and then decide whether $x\in K$.

To construct the automaton, let $\ell(h)$ denote the word length of
$h\in\Gr$ with respect to $\{a,b,c,d\}$. Since $a^2=b^2=c^2=d^2=bcd=1$, a shortest
word alternates the letter $a$ with letters in $\{b,c,d\}$, and taking sections we
obtain
\[
 \ell(h|_i)\le\frac{\ell(h)+1}{2}\qquad(i=0,1).
\]
Since $z\in K\le\St_\Gr(1)$, we have $x|_0=z|_0$ and $x|_1=z|_1x$. For an
automorphism of the form $hx$ with $h\in\Gr$, exactly one first-level section
retains the factor $x$, namely
\[
 T(h)x,\qquad T(h)=h|_{1^{h^{-1}}}\,z|_1,
\]
while the other section is $h|_i z|_0$, where $i^{h}=0$, and belongs to $\Gr$.
The length estimate gives
\[
 \ell(T(h))\le\frac{\ell(h)+1}{2}+\ell(z|_1).
\]
Consequently, all iterates of $T$ starting at $h=1$ have length at most
$2\ell(z|_1)+1$. We compute these iterates until one repeats, using the word
problem in $\Gr$ to detect equality. Only finitely many other sections are then
attached to this recursion, and each is a finite-state element of $\Gr$.
Together they give a finite automaton for $x$.

We next decide whether $x$ belongs to $\Gr$. Let $\mathcal N=\{1,a,b,c,d\}$,
which is closed under taking sections. Iterating the length estimate, we obtain
\[
 \ell(h|_v)\le 1+\frac{\ell(h)-1}{2^{|v|}}\qquad(h\in\Gr).
\]
Since word lengths are integers, every element of $\Gr$ has all its sufficiently
deep sections in $\mathcal N$. Using the decidability of
equality, let us mark the states of the automaton which are equal to elements of
$\mathcal N$. Since $\mathcal N$ is closed under sections, all sufficiently deep
sections belong to $\mathcal N$ if and only if the reachable unmarked states contain
no directed cycle, and this condition is decidable on the finite state graph. If
there is such a cycle, then $x\notin\Gr$. Otherwise, we decide membership in $\Gr$
for the reachable states by repeatedly choosing an unmarked state whose two
successor states have already been treated. Such a state exists because the
unmarked subgraph is acyclic. Suppose that its sections $h_0,h_1$ have been
found to belong to $\Gr$.
Multiplication on the right by $a\in\Gr$, if necessary, makes its root permutation
trivial without changing these sections. Since $\psi(\St_\Gr(1))\ge K\times K$, it belongs to $\Gr$ if and only if
\[
 (h_0K,h_1K)\in J:=\psi(\St_\Gr(1))/(K\times K)\le(\Gr/K)^2.
\]
The finite subgroup $J$ is computed from the sections of $b,c,d,b^a,c^a,d^a$, which
generate $\St_\Gr(1)$, and the cosets $h_iK$ are computed from the action of $h_i$
on $\Lv_{n_0}$. If a state turns out not to belong to $\Gr$, then $x\notin\Gr$ by
self-similarity.

Finally, if $x\in\Gr$, then $x\in K$ if and only
if the action of $x$ on $\Lv_{n_0}$ coincides with that of an element of $K$, which
can be checked in the finite quotient $\Gr/\St_\Gr(n_0)$.
\end{proof}

\section{Open questions}\label{sec:activity}

We conclude by discussing the hypotheses of our results and some open questions.
Bounded activity is used only in Lemma~\ref{lem:activity}, to rule out the equation
$x=\lambda(k)\rho(x)$. More generally, the proof works whenever, for some function
$f\colon\mathbb N\to(0,\infty)$, every $g\in G$ satisfies $A_n(g)=O(f(n))$ and some
$k\in K'$ satisfies $\sum_{j<n}A_j(k)\ne O(f(n))$. Theorem~\ref{thm:Cpoly} is the
case $f(n)=n^m$, and the long-range group (Theorem~\ref{thm:D}) shows that it
applies beyond bounded activity; see \cite{Sid00} for polynomial activity of
finite-state automorphisms and \cite{AAV13} for the amenability of groups generated
by finite-state automorphisms of linear activity.

Some activity assumption is however needed: the finitely
generated branch groups of \cite{SW03} contain non-abelian free subgroups, whereas
groups of bounded automorphisms, finite-state or not, do not \cite{Nek10} (those
generated by finite-state bounded automorphisms are even amenable \cite{BKN10});
hence the groups of \cite{SW03} admit no faithful action with bounded activity on a
regular rooted tree. Similarly, each generator of the Hanoi towers group $H^{(k)}$,
$k\ge4$, interchanges two first-level vertices and has itself as a section at the
other $k-2$ \cite{GS06}, so that its activity is $(k-2)^n$.

\begin{question}\label{ques:activity}
Is there a finitely generated branch group which is double coset separable? More
specifically, can the activity assumption or regular branching be omitted from
Theorem~\ref{thm:A}? Are the Hanoi towers groups $H^{(k)}$, $k\ge4$, double coset
separable? Is the second derived subgroup of the GGS-group with constant non-zero
defining vector finitely generated, so that Theorem~\ref{thm:C} applies to it
\cite{DFG23}?
\end{question}

When the branch kernel is non-trivial, as for $H^{(3)}$, or the group is not
branch, as for the long-range group, Theorem~\ref{thm:B} gives no information on the
closure of the product (Proposition~\ref{prop:hanoi},
Remark~\ref{rem:linear-witness}).

\begin{question}
Describe the closure of $PQ$ for $H^{(3)}$ and for the long-range group of
Theorem~\ref{thm:D}. Is the long-range group subgroup separable?
\end{question}

For the branching subgroup $K$ of Theorem~\ref{thm:A}, the twisted class $[1]_\rho$
is not closed, whereas the conjugacy classes of $\Gr$ are closed
\cite[p.~505]{Leo98}; by Proposition~\ref{prop:twisted}, closedness of
$D_\alpha(K)$ amounts to closedness of $\iota(\Gamma_\alpha)\iota(\Delta_K)$ in $G$.

\begin{question}
Is $D_\alpha(K)$ closed for every automorphism $\alpha$ of $K$? Is the first
Grigorchuk group twisted conjugacy separable with respect to all its automorphisms?
\end{question}

\begin{question}
Proposition~\ref{prop:membership} decides membership in the product $PQ$. Is
membership in arbitrary rational subsets of $\Gr$ decidable?
\end{question}
\subsection*{Acknowledgements}
Both authors are members of the Gruppo Nazionale per le Strutture Algebriche,
Geometriche e le loro Applicazioni (GNSAGA) of the Istituto Nazionale di Alta
Matematica (INdAM).

\subsection*{AI declaration}
OpenAI's ChatGPT was utilized for language support and computational checks of some selected identities and quotients during manuscript preparation. The authors take full responsibility for all mathematical content.

\end{document}